\documentclass[12pt]{article}

\usepackage{amsmath} 
\usepackage{amsthm,amsmath,mathrsfs,amsfonts,amssymb}
\usepackage{graphicx,psfrag,epsf,fixmath,caption,subcaption}
\usepackage{enumerate}
\usepackage{natbib} \RequirePackage[colorlinks,citecolor=blue,urlcolor=blue,linkcolor=blue]{hyperref}
\usepackage{url} 

\usepackage{enumitem} 

\newcommand{\blind}{1}
\newtheorem{assumption}{Assumption}

\newtheorem{theorem}{Theorem}

\theoremstyle{definition}

\begin{document}

\def\spacingset#1{\renewcommand{\baselinestretch}{#1}\small\normalsize} \spacingset{1}

\if1\blind
{
\title{\bf {\normalsize HETEROSKEDASTICITY-ROBUST INFERENCE IN LINEAR REGRESSION MODELS WITH MANY COVARIATES}}
\author{ {\small KOEN JOCHMANS}\thanks{Address: University of Cambridge, Faculty of Economics, Austin Robinson Building, Sidgwick Avenue, Cambridge CB3~9DD, United Kingdom. E-mail: \texttt{kj345@cam.ac.uk}. \newline An Associate Editor and three referees provided very constructive feedback on earlier versions of this paper and pointed out a non sequitur in the derivation of a primitive condition. Matias Cattaneo generously shared and discussed his replication material. I am most grateful to each of them for their help. Financial support from the European Research Council through grant n\textsuperscript{o} 715787 (MiMo) is gratefully acknowledged. \newline Revision history: September 5, 2018 (first version); August 28, 2019 (revision); April 27, 2020 (revision); \today  \ (this final version).}\\ {\small UNIVERSITY OF CAMBRIDGE}
}
\date{} 
  \maketitle
} \fi

\if0\blind
{
  \bigskip
  \bigskip
  \bigskip
  \begin{center}
    {\bf {\normalsize HETEROSKEDASTICITY-ROBUST INFERENCE IN LINEAR REGRESSION MODELS WITH MANY COVARIATES}}
\end{center}
  \medskip
} \fi

\medskip
\vspace{-.25cm}

\begin{abstract}
\noindent
We consider inference in linear regression models that is robust to heteroskedasticity and the presence of many control variables. When the number of control variables increases at the same rate as the sample size the usual heteroskedasticity-robust estimators of the covariance matrix are inconsistent. Hence, tests based on these estimators are size distorted even in large samples. An alternative covariance-matrix estimator for such a setting is presented that complements recent work by \cite{CattaneoJanssonNewey2018}. We provide high-level conditions for our approach to deliver (asymptotically) size-correct inference as well as more primitive conditions for three special cases. Simulation results and an empirical illustration to inference on the union premium are also provided.



\end{abstract}


\spacingset{1.45}

\renewcommand{\theequation}{\arabic{section}.\arabic{equation}} 
\setcounter{equation}{0}

\newpage

\section{Introduction}

When performing inference in linear regression models it is common practice to safeguard against (conditional) heteroskedasticity of unknown form. The estimator of the covariance matrix of the least-squares estimator proposed by \cite{Eicker1963,Eicker1967} and \cite{White1980} is known to be biased. The bias can be severe if the regressor design contains observations with high leverage \citep{ChesherJewitt1987}.\footnote{The leverage of an observation $i$ is defined as the $i$th diagonal element of the hat matrix, i.e., the matrix that transforms the observed outcomes into fitted values. It is bounded between zero and one and measures the influence of the observation on its own fitted value; a larger value reflects a higher influence (see, e.g., \citealt{HoaglinWelsch1978}).} A necessary condition for the least-squares estimator to be asymptotically normal is that maximal leverage vanishes as the sample size grows large \citep{Huber1973}. This condition, then, also implies consistency of the robust covariance-matrix estimator (under regularity conditions).

The requirement that maximal leverage vanishes is problematic when the regressors include a large set of control variables. 
Under asymptotics where their number, $q_n$, grows with the sample size, $n$, the robust covariance-matrix estimator will be inconsistent unless $q_n/n\rightarrow 0$, as demonstrated by \cite{CattaneoJanssonNewey2018}. 
They obtained the same result for the other members of the so-called HC-class of covariance-matrix estimators (see, e.g., \citealt{LongErvin2000} and \citealt{MacKinnon2012} for reviews) and showed that the jackknife variance estimator of \cite{MacKinnonWhite1985}, although inconsistent, can be used to perform asymptotically-conservative inference under asymptotics where $q_n/n\nrightarrow 0$. On the other hand, a bias-corrected covariance-matrix estimator in the spirit of \cite{HartleyRaoKiefer1969} and \cite{BeraSuprayitnoPremaratne2002} is shown to be consistent under conditions that bound maximal leverage below $\frac{1}{2}$. One implication of these conditions is that $\lim\sup_n q_n/n < \frac{1}{2}$. Although not guaranteed to be positive semi-definite (see \citealt{BeraSuprayitnoPremaratne2002} for a discussion), this estimator is attractive as it is asymptotically equivalent to a minimum-norm unbiased estimator in the sense of \cite{Rao1970}.

In this paper we discuss an alternative estimator of the covariance matrix that can deal with designs where maximal leverage is bounded away from 1. As such it remains consistent when $\lim\sup_n q_n/n < 1$. To show this we need to impose additional conditions relative to \cite{CattaneoJanssonNewey2018}. A consistency result is first provided under high-level conditions. Primitive conditions are then given for three special cases; the partially-linear regression model, the one-way model for short panel data, and the generic linear model with increasing dimension. Again, our covariance-matrix estimator need not be positive semi-definite in small samples. It is further not invariant to changes in the scale of the regression slopes. Achieving such invariance under asymptotics where $\lim\sup_n q_n/n > \frac{1}{2}$ appears difficult, as we discuss below.

The idea underlying our variance estimator can be traced back to work by  \citet[Remark 4 and Lemma 5]{KlineSaggioSoelvsten2018}; see below. The chief difference lies in the conditions under which the consistency result is obtained. They considered settings where the observations are independent and the regressors are fixed, and entertained models that are correctly specified with regression functions that are uniformly bounded. This is reasonable in analysis-of-variance problems, which are the main focus of application in their work. Here, we maintain the framework of \cite{CattaneoJanssonNewey2018}. Some dependence between observations is allowed, regressors are stochastic, and the model can feature (vanishing) mispecification bias. Our consistency result can be understood to be an extension of \cite{KlineSaggioSoelvsten2018} to such settings. It allows for regressors to have unbounded support and for observations to depend on a growing number of parameters. Our conditions do rule out dynamic models, for example. \cite{Verdier2018} provides estimation and inference results for two-way models estimated from short panel data that allow for dynamics over time.

In Section \ref{sec:1} we introduce the framework, present our covariance-matrix estimator, and provide a consistency result under a set of high-level conditions. In Section \ref{sec:2} we connect our work to the literature and, notably, to \cite{CattaneoJanssonNewey2018} and \cite{KlineSaggioSoelvsten2018}. In Section \ref{sec:3} we provide primitive conditions for three special cases of our setup. In Section \ref{sec:4} we present and discuss the results from Monte Carlo experiments and apply our variance estimator to perform inference on the union wage premium. A short conclusion ends the paper. The supplemental appendix contains technical details and additional simulation results.

\section{Inference with many regressors} \label{sec:1}

\subsection{Framework}

Consider the linear model
\begin{equation} \label{eq:model}
y_{i,n} = \boldsymbol{x}_{i,n}^\prime\boldsymbol{\beta} + \boldsymbol{w}_{i,n}^\prime \boldsymbol{\gamma}_n + u_{i,n}, \qquad i=1,\ldots, n,
\end{equation}
where $y_{i,n}$ is a scalar outcome, $\boldsymbol{x}_{i,n}$ is a vector of regressors of fixed dimension $r$, $\boldsymbol{w}_{i,n}$ is a vector of covariates whose dimension, $q_n$, may grow with $n$, and $u_{i,n}$ is an unobserved error term. Our aim is to perform asymptotically-valid inference on $\boldsymbol{\beta}$ that is robust to (conditional) heteroskedasticity, when $\boldsymbol{\gamma}_n$ is high-dimensional, in the sense that $q_n$ is not a vanishing fraction of the sample size. In such a case, the nuisance parameter $\boldsymbol{\gamma}_n$ is not consistently estimable, in general.


The (ordinary) least-squares estimator of $\boldsymbol{\beta}$ is
$$
\boldsymbol{\hat\beta}_n := 
\left(\sum_{i=1}^n \boldsymbol{\hat{v}}_{i,n}\boldsymbol{\hat{v}}_{i,n}^\prime\right)^{-1}
\left(\sum_{i=1}^n \boldsymbol{\hat{v}}_{i,n} y_{i,n}\right),
$$
where
$$
\boldsymbol{\hat v}_{i,n} := \sum_{j=1}^n (\boldsymbol{M}_n)_{i,j} \, \boldsymbol{x}_{j,n},
\qquad
(\boldsymbol{M}_n)_{i,j}:= \lbrace i = j\rbrace - 
\boldsymbol{w}_{i,n}^\prime \left(\sum_{k=1}^n\boldsymbol{w}_{k,n} \boldsymbol{w}_{k,n}^\prime \right)^{-1} \boldsymbol{w}_{j,n},
$$
and $\lbrace \cdot \rbrace$ denotes the indicator function. We will provide an inference result based on the limit distribution of $\boldsymbol{\hat \beta}_n$.
We begin by stating a set of high-level conditions that guarantee this distribution to be Gaussian that cover the case where $q_n/n \nrightarrow 0$ as $n\rightarrow\infty$. Our point of departure for this is \citet[Theorem 1]{CattaneoJanssonNewey2018} and we closely follow their notation.

Let $\mathcal{X}_n:=(\boldsymbol{x}_{1,n},\ldots,\boldsymbol{x}_{n,n})$ and let $\mathcal{W}_n$ denote a collection of random variables such that $\mathbb{E}(\boldsymbol{w}_{i,n} \vert \mathcal{W}_n) = \boldsymbol{w}_{i,n}$. We introduce 
$$
\varepsilon_{i,n}:= u_{i,n} - e_{i,n},
\qquad
\boldsymbol{V}_{i,n}:= \boldsymbol{x}_{i,n} - \mathbb{E}(\boldsymbol{x}_{i,n} \vert \mathcal{W}_n),
$$
where $e_{i,n}:=\mathbb{E}(u_{i,n} \vert \mathcal{X}_n, \mathcal{W}_n)$, to state our first assumption. We use $\lvert \cdot \rvert$ to denote the cardinality of a set.

\begin{assumption}[Sampling] \label{ass:sampling}
The errors $\varepsilon_{i,n}$ are uncorrelated across $i$ conditional on $\mathcal{X}_n$ and $\mathcal{W}_n$, and the collections $\lbrace \varepsilon_{i,n}, \boldsymbol{V}_{i,n} : i\in N_g \rbrace$ are independent across $g$ conditional on $\mathcal{W}_n$, where $\lbrace N_1,\ldots, N_{G_n} \rbrace$ represents a partition of $\lbrace 1,\ldots, n \rbrace$ into $G_n$ sets such that $\max_g \lvert N_g \rvert = O(1)$.
\end{assumption}

\noindent
This assumption covers standard randomly-sampled data but also repeated-measurement data (such as short panel data) where strata are independent but dependence between observations within the strata is allowed, for example. 

The second assumption contains regularity conditions. We let 
$$
\sigma_{i,n}^2:=\mathbb{E}(\varepsilon_{i,n}^2 \vert \mathcal{X}_n, \mathcal{W}_n),
\qquad
\boldsymbol{\tilde V}_{i,n}:= \sum_{j=1}^n (\boldsymbol{M}_n)_{i,j} \, \boldsymbol{V}_{j,n},
$$
denote the Euclidean and Frobenius norms by $\lVert \cdot \rVert$, and write $\lambda_{\min}(\cdot)$ for the minimum eigenvalue of its argument. 

\begin{assumption}[Design] \label{ass:design}
With probability approaching one $\sum_{i=1}^n \boldsymbol{w}_{i,n}\boldsymbol{w}_{i,n}^\prime$ has full rank,
$$
\max_i 
\left(\mathbb{E}(\varepsilon_{i,n}^4\vert \mathcal{X}_n,\mathcal{W}_n)
+ \frac{1}{\sigma_{i,n}^2}+\mathbb{E}(\lVert \boldsymbol{V}_{i,n}  \rVert^4 \, \vert \mathcal{W}_n ) \right) + \frac{1}{\lambda_{\min}\left( \frac{\sum_{i=1}^n \mathbb{E}(\boldsymbol{\tilde V}_{i,n} \boldsymbol{\tilde V}_{i,n}^\prime \vert \mathcal{W}_n)}{n}\right)}
=
O_p(1),
$$ 
and $\lim\sup_n q_n/n <1$.
\end{assumption}

\noindent
The rank condition on the design matrix $\sum_{i=1}^n \boldsymbol{w}_{i,n}\boldsymbol{w}_{i,n}^\prime$ is standard. Furthermore, given that the slope coefficients on $\boldsymbol{w}_{i,n}$ are not of direct interest to us, dropping any covariates that are (perfectly) collinear is not an issue. The second condition contains conventional moment conditions. The third condition, finally, allows for $q_n$ to grow at the same rate as the sample size.

Our setting covers situations where the regression in \eqref{eq:model} is a linear-in-parameters mean-square approximation to the conditional expectation $\mu_{i,n}:=\mathbb{E}(y_{i,n} \vert \mathcal{X}_n, \mathcal{W}_n)$, in the sense that we allow that $e_{i,n}\neq 0$. The third assumption contains conditions on how fast such an approximation should improve. They are expressed in terms of the two constants
$$
\varrho_n
 := 
\frac{\sum_{i=1}^n \mathbb{E}(e_{i,n}^2)}{n} ,
\qquad
\rho_{n}
 := 
\frac{\sum_{i=1}^n \mathbb{E}(\mathbb{E}(e_{i,n} \vert \mathcal{W}_n)^2)}{n}.
$$
The assumption also contains a similar restriction on how well $\boldsymbol{V}_{i,n}$ can be approximated by
$$
\boldsymbol{v}_{i,n}  := \boldsymbol{x}_{i,n} - 
\left(\sum_{j=1}^n \mathbb{E}(\boldsymbol{x}_{j,n}\boldsymbol{w}_{j,n}^\prime )\right) 
 \left(\sum_{j=1}^n \mathbb{E}(\boldsymbol{w}_{j,n} \boldsymbol{w}_{j,n}^\prime)\right)^{-1}
\boldsymbol{w}_{i,n},
$$
the deviation of $\boldsymbol{x}_{i,n}$ from its population linear projection. This is expressed using the constant
$$
\chi_n := \frac{\sum_{i=1}^n \mathbb{E}(\lVert\boldsymbol{Q}_{i,n}\rVert^2)}{n},
$$
where
$
\boldsymbol{Q}_{i,n}:= \mathbb{E}(\boldsymbol{v}_{i,n} \vert \mathcal{W}_n).
$

\begin{assumption}[Approximations] \label{ass:approximations}
$\chi_n = O(1)$, $\varrho_n + n(\varrho_n-\rho_n) + n\chi_n \varrho_n = o(1)$, and $\max_i \, \lVert  \boldsymbol{\hat{v}}_{i,n} \rVert/\sqrt{n} = o_p(1)$.
\end{assumption}	

\noindent
The last part of this assumption is a high-level negligibility condition on the residuals from the auxiliary regression of $\boldsymbol{x}_{i,n}$ on $\boldsymbol{w}_{i,n}$. Given that
$$
\max_i
\left(\boldsymbol{\hat v}_{i,n}^\prime \left(\sum_{j=1}^n \boldsymbol{\hat v}_{j,n} \boldsymbol{\hat v}_{j,n}^\prime\right)^{-1} \, \boldsymbol{\hat v}_{i,n}\right)
\leq
\max_i \frac{\lVert \boldsymbol{\hat v}_{i,n} \rVert^2}{n} \, 
\left\lVert \left(\frac{\sum_{j=1}^n \boldsymbol{\hat v}_{j,n} \boldsymbol{\hat v}_{j,n}^\prime}{n}\right)^{-1} \right\rVert^2
$$
and should vanish in large samples for $\boldsymbol{\hat\beta}_n$ to be asymptotically normal \citep{Huber1973}, this requirement appears close to minimal.

Assumptions \ref{ass:sampling}--\ref{ass:approximations} coincide with Assumptions 1--3 in \cite{CattaneoJanssonNewey2018}. Consequently, by their Theorem 1,
\begin{equation} \label{eq:normality}
\boldsymbol{\Omega}_n^{-1/2} (\boldsymbol{\hat\beta}_n-\boldsymbol{\beta})
\overset{d}{\rightarrow}
\boldsymbol{N}(\boldsymbol{0},\boldsymbol{I}_{r})
\end{equation}
as $n\rightarrow\infty$, where
$$
\boldsymbol{\Omega}_n := 
\left(\sum_{i=1}^n{\boldsymbol{\hat{v}}}_{i,n}{\boldsymbol{\hat{v}}}_{i,n}^\prime\right)^{-1}
\left(\sum_{i=1}^n{\boldsymbol{\hat{v}}}_{i,n}{\boldsymbol{\hat{v}}}_{i,n}^\prime \, \sigma_{i,n}^2\right)
\left(\sum_{i=1}^n{\boldsymbol{\hat{v}}}_{i,n}{\boldsymbol{\hat{v}}}_{i,n}^\prime\right)^{-1},
$$
and $\boldsymbol{I}_r$ denotes the $r\times r$ identity matrix.

\subsection{Variance estimation}

Constructing confidence intervals and test statistics based on \eqref{eq:normality} requires an estimator of $\boldsymbol{\Omega}_n$, and thus of
$$
\boldsymbol{\Sigma}_n := 
\sum_{i=1}^n{\boldsymbol{\hat{v}}}_{i,n}{\boldsymbol{\hat{v}}}_{i,n}^\prime \, \sigma_{i,n}^2.
$$
When $q_n/n \nrightarrow 0$ and the errors are permitted to be heteroskedastic, the construction of a consistent estimator is non-trivial. To appreciate the problem, consider the estimator of \cite{Eicker1963,Eicker1967} and \cite{White1980}, which uses  
$$
\boldsymbol{\hat\Sigma}_n := 
\sum_{i=1}^n{\boldsymbol{\hat{v}}}_{i,n}{\boldsymbol{\hat{v}}}_{i,n}^\prime \, {\hat u}_{i,n}^2,
$$
where 
$
\hat{u}_{i,n}:= \sum_{j=1}^n (\boldsymbol{M}_{n})_{i,j} \, 
(y_{j,n} - \boldsymbol{x}_{j,n}^\prime \boldsymbol{\hat \beta}_n)
$
are the least-squares residuals. This estimator is well known to be (conditionally) biased. The bias arises from the sampling noise in the least-squares estimator and can be severe \citep{ChesherJewitt1987}. 
Unless $q_n/n \rightarrow 0$ as $n\rightarrow \infty$, some observations will remain influential, in the sense that maximal leverage does not vanish. This causes the bias in $\boldsymbol{\hat{\Sigma}}_n$ to persist in large samples, implying that it is inconsistent.

The alternative to $\boldsymbol{\hat\Sigma}_n$ that we consider in this paper is
$$
\boldsymbol{\acute \Sigma}_n:= 
\sum_{i=1}^n{\boldsymbol{\hat{v}}}_{i,n}{\boldsymbol{\hat{v}}}_{i,n}^\prime \, (y_{i,n}\acute{u}_{i,n}) ,
\qquad
\acute{u}_{i,n} := \frac{\hat{u}_{i,n}}{(\boldsymbol{M}_n)_{i,i}}.
$$
As stated, this estimator is well defined provided that 
$$
\min_i (\boldsymbol{M}_n)_{i,i} > 0.
$$
Notice that $(\boldsymbol{M}_n)_{i,i}=0$ means that the model reserves a parameter for this observation. This implies that the auxiliary regression of the regressors of interest on the other covariates yields a perfect prediction, in the sense that $\boldsymbol{\hat v}_{i,n}= \boldsymbol{0}$. Consequently, such an observation does not carry information on $\boldsymbol{\beta}$ and can be dropped. It does not affect the least-squares estimator $\boldsymbol{\hat\beta}_n$ and does not contribute to its covariance matrix $\boldsymbol{\Omega}_n$.
This is important as perfect prediction of this form arises frequently in empirical work when many dummy variables are included.


Additional conditions are needed to show that $\boldsymbol{\acute{\Sigma}}_n$ is consistent. We let
$$
\boldsymbol{\tilde Q}_{i,n} := \sum_{j=1}^n (\boldsymbol{M}_n)_{i,j} \, \boldsymbol{Q}_{i,n}.
$$
and collect one such set of conditions in the following assumption.

\begin{assumption}[Variance estimation] \label{ass:inference}
$n \varrho_n = O(1)$, $\Pr(\min_i (\boldsymbol{M}_n)_{i,i} > 0 )\rightarrow 1$, 
$$
\frac{1}{\min_i (\boldsymbol{M}_n)_{i,i}} = O_p(1),
\qquad
\frac{\sum_{i=1}^n \lVert \boldsymbol{\tilde Q}_{i,n} \rVert^4 }{n}
= O_p(1),
$$
and $\max_i \lVert \mu_{i,n} \rVert /\sqrt{n} = o_p(1)$.
\end{assumption}

\noindent
The first part of Assumption \ref{ass:inference} is a small-bias condition; it is relevant only when \eqref{eq:model} is misspecified, in the sense that $e_{i,n}\neq 0$. In that case it is a strengthening of Assumption \ref{ass:approximations} only when $\chi_n = o(1)$. The conditions on the diagonal entries of the projection matrix are very weak. Providing primitive conditions for them in great generality appears to be difficult. However, when $\boldsymbol{w}_{i,n}$ is multivariate normal they follow under $\lim\sup_n q_n/n<1$ as stated in Assumption \ref{ass:design} in the same way as in \cite{CattaneoJanssonNewey2018}. In the one-way panel model they hold automatically while \cite{Verdier2018} gives sufficient conditions for them to be satisfied in the two-way model. To understand why the last part of Assumption \ref{ass:inference} is needed, note that the (conditional) variance of $y_{i,n}\, \acute{u}_{i,n}$ depends on $\mu_{i,n}^2$. The requirement that $\max_i \mu_{i,n}^2 = o_p(n)$ allows to control the variance of 
$
\boldsymbol{\acute \Sigma}_n.
$
Weak moment requirements typically suffice for this condition to be satisfied. The condition on $\boldsymbol{\tilde Q}_{i,n}$ is used in concordance with the condition on $\mu_{i,n}$. One simple sufficient condition for it is that $n \chi_n = O(1)$, but it can also be satisfied when $\chi_n = O(1)$. Primitive conditions for Assumption \ref{ass:inference} in three special cases are given below.

We can now state our consistency result.

\begin{theorem}[Inference] \label{thm:inference}
Let Assumptions \ref{ass:sampling}--\ref{ass:inference} hold. Then
$$
\boldsymbol{\Sigma}_n^{-1}\boldsymbol{\acute\Sigma}_n \overset{p}{\rightarrow} \boldsymbol{I}_r
$$	
as $n\rightarrow\infty$.
\end{theorem}	

\noindent
Theorem \ref{thm:inference}, combined with the limit result in \eqref{eq:normality}, implies that
$$
\boldsymbol{\acute\Omega}_n^{-1/2} (\boldsymbol{\hat\beta}_n-\boldsymbol{\beta})
\overset{d}{\rightarrow}
\boldsymbol{N}(\boldsymbol{0},\boldsymbol{I}_{r})
$$
as $n\rightarrow\infty$, where
$$
\boldsymbol{\acute\Omega}_n := 
\left(\sum_{i=1}^n{\boldsymbol{\hat{v}}}_{i,n}{\boldsymbol{\hat{v}}}_{i,n}^\prime\right)^{-1}
\left(\sum_{i=1}^n{\boldsymbol{\hat{v}}}_{i,n}{\boldsymbol{\hat{v}}}_{i,n}^\prime \, (y_{i,n}\, \acute{u}_{i,n})\right)
\left(\sum_{i=1}^n{\boldsymbol{\hat{v}}}_{i,n}{\boldsymbol{\hat{v}}}_{i,n}^\prime\right)^{-1}.
$$
This result permits the construction of test statistics that (in large samples) will have correct size and of confidence regions that will exhibit correct coverage.

\section{Connections to the literature} \label{sec:2}

\paragraph{HC-class estimators.}
The bias in the \cite{Eicker1963,Eicker1967} and \cite{White1980} estimator has led to a variety of modifications to it being proposed that, following \cite{MacKinnonWhite1985}, are often referred to as the HC-class of covariance-matrix estimators. These estimators are reviewed in \cite{LongErvin2000} and \cite{MacKinnon2012}. Unfortunately, as shown by \citet[Theorem 3]{CattaneoJanssonNewey2018}, none of these alternatives is consistent, in general, when $q_n/n\nrightarrow 0$. We briefly review their main findings on these estimators here.

The first variance estimator, HC1, differs from the conventional estimator, HC0, in that it performs a degrees-of-freedom correction \citep[Eq.~2.11]{Hinkley1977}. This estimator will be consistent in the special case where errors are homoskedastic and the covariate design is balanced, i.e., when $(\boldsymbol{M}_n)_{1,1}=\ldots=(\boldsymbol{M}_n)_{n,n}$.

The second variance estimator, HC2, uses\footnote{We follow \cite{CattaneoJanssonNewey2018} and construct the HC-class variance estimators using the projection matrix $\boldsymbol{M}_n$. The original proposals were made in a context where $q_n$ is treated as fixed and did not differentiate between $\boldsymbol{x}_{i,n}$ and $\boldsymbol{w}_{i,n}$; they used the annihilator matrix that projects out both sets of variables. The difference is asymptotically negligible under our assumptions.
}
$$
\sum_{i=1}^n \boldsymbol{\hat{v}}_{i,n} \boldsymbol{\hat{v}}_{i,n}^\prime 
(\hat{u}_{i,n}\, \acute{u}_{i,n})
$$
as an estimator of $\boldsymbol{\Sigma}_n$ \citep{HornHornDuncan1975}. This estimator will be consistent under homoskedasticity. Because
$$
y_{i,n} \acute{u}_{i,n} 
=
\hat{u}_{i,n}\, \acute{u}_{i,n}
+
\hat{y}_{i,n}\, \acute{u}_{i,n},
$$
where $\hat{y}_{i,n}:=y_{i,n}-\hat{u}_{i,n}$ are fitted values, 
$\boldsymbol{\acute{\Sigma}}_n$ can be interpreted as a bias-corrected version of the HC2 estimator.

The third variance estimator, HC3, is constructed with
$$
\sum_{i=1}^n \boldsymbol{\hat{v}}_{i,n} \boldsymbol{\hat{v}}_{i,n}^\prime
(\acute{u}_{i,n}\, \acute{u}_{i,n}) 
$$
(\citealt{MacKinnonWhite1985}). While this estimator is inconsistent, its probability limit exceeds $\boldsymbol{\Sigma}_n$ (in the matrix sense). It follows that (under Assumptions \ref{ass:sampling}--\ref{ass:approximations}) test procedures based on HC3 will be asymptotically conservative when $\lim\sup_n q/n \in (0,1)$, both under homoskedasticity and heteroskedasticity. 

\paragraph{Bias-corrected estimation.}
\cite{CattaneoJanssonNewey2018} also considered the estimator
$$
\boldsymbol{\grave \Sigma}_n:= 
\sum_{i=1}^n{\boldsymbol{\hat{v}}}_{i,n}{\boldsymbol{\hat{v}}}_{i,n}^\prime \, \left(\sum_{j=1}^n ((\boldsymbol{M}_n * \boldsymbol{M}_n)^{-1})_{i,j} \, {\hat u}_{j,n}^2\right),
$$
where $\boldsymbol{M}_n * \boldsymbol{M}_n$ denotes the elementwise product of the matrix $\boldsymbol{M}_n$. This estimator has its origins in work by 
\cite{HartleyRaoKiefer1969} and \cite{Rao1970} and can be motivated through an (asymptotic) bias calculation; see also \cite{BeraSuprayitnoPremaratne2002}, and \cite{Anatolyev2018} for a refinement under homoskedasticity. For $\boldsymbol{\grave \Sigma}_n$ to be well defined, $\boldsymbol{M}_n * \boldsymbol{M}_n$ needs to be nonsingular. Necessary and sufficient conditions for this to be the case are stated in \cite{Mallela1972} but these are neither simple nor intuitive \citep{HornHornDuncan1975}. As noted by \cite{HornHorn1975}, a simple sufficient condition is that
$$
\min_i (\boldsymbol{M}_n)_{i,i} > \frac{1}{2}.
$$
Depending on the problem at hand it may also be necessary; an example is the one-way panel model.

\citet[Theorem 4]{CattaneoJanssonNewey2018} show that $\boldsymbol{\grave \Sigma}_n$ is consistent for $\boldsymbol{\Sigma}_n$ if
\begin{equation} \label{eq:largeM}
\Pr\left(\min_i (\boldsymbol{M}_n)_{i,i} > \frac{1}{2}\right)\rightarrow 1,
\qquad
\frac{1}{\min_i (\boldsymbol{M}_n)_{i,i}-\frac{1}{2}} = O_p(1),
\end{equation}
are added to Assumptions \ref{ass:sampling}--\ref{ass:approximations}. 
Because $\sum_{i=1}^n (\boldsymbol{M}_n)_{i,i} = n-q_n$, 
$\min_i (\boldsymbol{M}_n)_{i,i}\leq 1-q_n/n$, and so
$$
\lim\sup_n q_n/n <\frac{1}{2}
$$
is required for \eqref{eq:largeM} to be satisfied. This, in turn, is a strengthening of the condition that $\lim\sup_n q_n/n < 1$ in Assumption \ref{ass:design}.

\cite{StockWatson2008} proposed a covariance-matrix estimator for linear fixed-effect models that is applicable to short panel data. It is based on an explicit calculation of the probability limit of $\boldsymbol{\hat{\Sigma}}_n-\boldsymbol{\Sigma}_n$, adjusting the inconsistent $\boldsymbol{\hat{\Sigma}}_n$ by subtracting from it a plug-in estimator of this limit quantity.
As discussed in \cite{CattaneoJanssonNewey2018}, $\boldsymbol{\grave \Sigma}_n$ can be understood to be a generalization of this approach to the generic setting where $q_n/n \nrightarrow 0$.

Like $\boldsymbol{\acute{\Sigma}}_n$, the estimator $\boldsymbol{\grave \Sigma}_n$ need not be positive semi-definite in small samples. Being based on estimators of the individual $\sigma_{i,n}^2$ that are linear combinations of $\hat{u}_{1,n}^2,\ldots,\hat{u}_{n,n}^2$ it does, however, retain the invariance to changes in $\mu_{i,n}$ that is inherent in the HC-class estimators. This is a desirable feature and explains why a restriction on the magnitude of $\mu_{i,n}$ is not needed for this estimator to be consistent. Furthermore, jointly estimating $\sigma_{1,n}^2,\ldots,\sigma^2_{n,n}$ in this way is attractive from an efficiency point of view, as it is asymptotically equivalent to a minimum-norm unbiased estimator (see \citealt{HartleyRaoKiefer1969} and \citealt{Rao1970} for details). 


\paragraph{Inference on variance components.}
The variance estimator $\boldsymbol{\acute{\Sigma}}_n$ is closely related to the work of \cite{KlineSaggioSoelvsten2018}. In our context, their proposal is to estimate each $\sigma_{i,n}^2$ by the cross-fit \citep{NeweyRobins2018} estimator
$$
y_{i,n} \, \check{u}_{i,n},
$$
where $\check{u}_{i,n}$ is the residual for observation $i$ when the slope coefficients are estimated from the sample from which the $i$th observation has been omitted. When the regression model is correctly specified---i.e., when $e_{i,n}=0$---then $y_{i,n} \,\check{u}_{i,n}$ is (conditionally) unbiased, provided that the leave-own-out regression slopes are well defined.
This will be the case if maximal leverage is bounded away from one.\footnote{This is a slight strengthening of our requirement in Assumption \ref{ass:inference} as we project out only $\boldsymbol{w}_{i,n}$ to obtain $\acute{u}_{i,n}$ while \cite{KlineSaggioSoelvsten2018} project out both the regressors $\boldsymbol{x}_{i,n}$ and the covariates $\boldsymbol{w}_{i,n}$ to obtain $\check{u}_{i,n}$.} In this case, from the Sherman-Morrison formula (as in, e.g., \citealt{Miller1974}),  
$$
y_{i,n} \, \check{u}_{i,n}
=
y_{i,n} \acute{u}_{i,n} + o_p(1),
$$
under our assumptions, connecting the cross-fit estimator of \cite{KlineSaggioSoelvsten2018} to our approach.

\cite{KlineSaggioSoelvsten2018} use this device to construct unbiased estimators of quadratic forms and present several applications. One of these (\citealt[Remark 4 and Lemma 5]{KlineSaggioSoelvsten2018}) is a consistency result for the implied covariance-matrix estimator 
$$
\sum_{i=1}^n \boldsymbol{\hat{v}}_{i,n}\boldsymbol{\hat{v}}_{i,n}^\prime \,
(y_{i,n} \,  \check{u}_{i,n})
$$
in designs with fixed regressors where errors are independent and maximal leverage is bounded away from unity. This result is established under the assumption that $e_{i,n}=0$ and that $\max_i \lVert \mu _{i,n}\rVert^2 = O(1)$, together with standard regularity conditions as those in Assumption \ref{ass:design}. It follows that our variance estimator can be seen as a modification of theirs that is targeted to a setting with many control variables. Theorem \ref{thm:inference} may then be understood to be an extension of their Lemma 5 to settings with stochastic regressors and (vanishing) mispecification bias, where the regressors can have unbounded support and observations may depend on a growing number of parameters.

The implied interpretation of $\boldsymbol{\acute\Sigma}_n$ as an (approximate) cross-fit estimator is also useful to highlight an apparent tension between invariance of a covariance-matrix estimator to changes in $\mu_{i,n}$ and the possibility for it to be consistent when we have $\lim\sup_n q_n/n>\frac{1}{2}$. A cross-fit estimator of $\sigma_{i,n}^2$ that is invariant is necessarily of the form $\dot{u}_{i,n} \ddot{u}_{i,n}$, where $\dot{u}_{i,n}$ and $\ddot{u}_{i,n}$ are two  least-squares residuals, each one obtained from a (conditionally) independent subsample that excludes the $i$th observation. Because these auxiliary regressions can be based on at most $\lfloor (n-1)/2 \rfloor$ and $\lceil (n-1)/2 \rceil$ observations such an approach cannot accommodate situations where $q_n/n>\frac{1}{2}$. The alternative estimators $y_{i,n} \, \check{u}_{i,n}$ and $y_{i,n} \, \acute{u}_{i,n}$ circumvent the need for two independent estimators of $u_{i,n}$ by using the level of the outcome variable, $y_{i,n}$, as a proxy for $u_{i,n}$. This allows to deal with cases where $\lim\sup_n q_n/n > \frac{1}{2}$ but makes the variance estimator sensitive to $\mu_{i,n}$.

\section{Examples} \label{sec:3}

We next provide more primitive conditions in three special cases that fit out general setup. We focus on sufficient conditions for Assumption \ref{ass:inference}. \cite{CattaneoJanssonNewey2018} already gave such conditions for the other assumptions---and, notably, for Assumption \ref{ass:approximations}---to hold.

\paragraph{Partially-linear model}
Suppose that observations on $(y_i, \boldsymbol{x}_i, \boldsymbol{z}_i)$ are independent and identically distributed.
The partially-linear regression model states that
$$
y_i = \boldsymbol{x}_i^\prime  \boldsymbol{\beta} + \varphi(\boldsymbol{z}_i) + \varepsilon_{i},
\qquad
\mathbb{E}(\varepsilon_i \vert \boldsymbol{x},\boldsymbol{z}_i) = 0, 
$$
for an unknown function $\varphi$. 
A series approximation of order $\kappa_n$ of $\varphi(\boldsymbol{z}_i)$ takes the form $\boldsymbol{w}_{i,n}^\prime\boldsymbol{\gamma}_n$ for $\boldsymbol{w}_{i,n}=(p_1(\boldsymbol{z}_i),\ldots, p_{\kappa_n}(\boldsymbol{z}_i))^\prime$ and $\lbrace p_1,\ldots, p_{\kappa_n} \rbrace$ a collection of basis functions such as orthogonal polynomials. Our estimator $\boldsymbol{\hat\beta}_n$, then, is the least-squares estimator of $\boldsymbol{\beta}$ in
$$
y_{i,n} = \boldsymbol{x}_{i}^\prime \boldsymbol{\beta} + \boldsymbol{w}_{i,n}^\prime\boldsymbol{\gamma}_n + u_{i,n},
\qquad
u_{i,n} = \varepsilon_{i} + \varphi(\boldsymbol{z}_i)-\boldsymbol{w}_{i,n}^\prime\boldsymbol{\gamma}_n.
$$
Note that $\mathbb{E}(u_{i,n} \vert \boldsymbol{x}_{i},\boldsymbol{z}_i) \neq 0$, in general. Consistency of $\boldsymbol{\hat\beta}_n$ requires that $\kappa_n\rightarrow \infty$ (and, thus, that $q_n\rightarrow\infty$) as $n\rightarrow\infty$.
Here, 
$$
\varrho_n = 
\min_{\boldsymbol{\gamma}}
\mathbb{E}\left(
\lVert
\varphi(\boldsymbol{z}_i)
-
\boldsymbol{w}_{i,n}^\prime \boldsymbol{\gamma}
\rVert^2\right), 
\qquad
\chi_n = 
\min_{\boldsymbol{\delta}}
\mathbb{E}\left(
\lVert
\mathbb{E}(\boldsymbol{x}_{i,n} \vert \boldsymbol{w}_{i,n})
-
\boldsymbol{\delta}^\prime \boldsymbol{w}_{i,n}
\rVert^2\right).
$$
In this example the number of covariates can be large when the dimension of $\boldsymbol{z}_i$ is large, so that many terms are included in the approximation even for small $\kappa_n$, or when the underlying functions are not (assumed to be) very smooth, so that a large $\kappa_n$ needs to be used to control bias.

Standard smoothness conditions on $\varphi$ imply that $n\varrho_n = O(1)$  \citep{Newey1997}, yielding the first condition of Assumption \ref{ass:inference}. The fourth condition can be validated in the same way, by imposing sufficient smoothness on $\mathbb{E}(\boldsymbol{x}_{i,n} \vert \boldsymbol{w}_{i,n})$ to get $n\chi_n = O(1)$, which implies that $\sum_{i=1}^n \lVert \boldsymbol{\tilde{Q}}_{i,n} \rVert^4 = O_p(n)$. This is different from (and stronger than) the first set of primitive conditions discussed by \cite{CattaneoJanssonNewey2018} to validate Assumption \ref{ass:approximations}, where more smoothness in $\mathbb{E}(\boldsymbol{x}_{i,n} \vert \boldsymbol{w}_{i,n})$ can be used to compensate for less smoothness in $\varphi$. Alternatively, if $\chi_n \leq \mathbb{E}(\lVert\mathbb{E}(\boldsymbol{x}_{i,n} \vert \boldsymbol{w}_{i,n})\rVert^2)=O(1)$ and a partitioning estimator (\citealt{CattaneoFarrell2013}) is used to approximate $\varphi$ then $\boldsymbol{M}_n$ is a band matrix. This is also sufficient to reach the desired result. Moreover, in this case, the first and fourth condition of Assumption \ref{ass:inference} are implied by the rate requirements in Assumption \ref{ass:approximations}, and by the second set of restrictions for this example given in \cite{CattaneoJanssonNewey2018}. Next, simple sufficient conditions for the requirement that $\max_i \lVert \mu_{i,n}\rVert/\sqrt{n} = o_p(1)$ are the moment conditions $\mathbb{E}(\lVert\boldsymbol{x}_i \rVert^{2+\theta})=O(1)$ and $\mathbb{E}(\lVert\varphi(\boldsymbol{z}_i) \rVert^{2+\theta})=O(1)$ for some $\theta>0$. While these two conditions are not imposed in \cite{CattaneoJanssonNewey2018}, they would not appear to be overly strong. 

\paragraph{One-way model for panel data}
For double-indexed data $(y_{(g,m)}, \boldsymbol{x}_{(g,m)})$, the fixed-effect model is
$$
y_{(g,m)} = \boldsymbol{x}_{(g,m)}^\prime \, \boldsymbol{\beta} + \alpha_g + \varepsilon_{(g,m)},
\qquad
g=1,\ldots, G_n, \quad  m=1\ldots,M,
$$
where $\alpha_g$ is a group-specific intercept and we assume that $\mathbb{E}(\varepsilon_{(g,m)} \vert \boldsymbol{x}_{(g,1)},\ldots, \boldsymbol{x}_{(g,M)}) = 0$. The regressors $\boldsymbol{x}_{(g,m)}$ are assumed independent between groups but may be dependent within each group. The errors $\varepsilon_{(g,m)}$ are independent between groups and (conditionally) uncorrelated within groups. The usual asymptotic embedding here has $G_n\rightarrow \infty$ with $M=O(1)$. The number of fixed effects grows at the same rate as the sample size; we have $n=G_n\times M$ and $q_n = G_n$ so that $q_n/n=1/M$, which does not vanish. These conditions fit Assumption \ref{ass:sampling}. 

The fixed-effect estimator is the least-squares estimator of $y_{(g,m)}$ on $\boldsymbol{x}_{(g,m)}$ and $G_n$ dummy variables that indicate group membership. The estimated coefficients on these dummies are computed from $M$ observations and are not consistent under our asymptotic approximation. In this example $\boldsymbol{M}_n = \boldsymbol{I}_{G_n} \otimes \boldsymbol{T}_M$, where 
$
(\boldsymbol{T}_M)_{m,m^\prime}
:=
\lbrace m = m^\prime \rbrace - M^{-1}
$
is the $M\times M$ matrix that transforms observations into deviations from their within-group mean. Consequently, 
$$
\boldsymbol{\acute{\Sigma}}_n
=
\frac{M}{M-1}
\sum_{g=1}^{G_n} \sum_{m=1}^M 
\boldsymbol{\tilde{x}}_{(g,m)} \, \boldsymbol{\tilde{x}}_{(g,m)}^\prime \ 
y_{(g,m)} \, (\tilde{y}_{(g,m)}-\boldsymbol{\tilde{x}}_{(g,m)}^\prime \boldsymbol{\hat{\beta}}_n),
$$
where
$
\tilde{y}_{(g,m)} : = \sum_{m^\prime = 1}^{M} (\boldsymbol{T}_M)_{m,m^\prime} \, y_{(g,m^\prime)}
$,
and $\boldsymbol{\tilde{x}}_{(g,m)}$ and $\tilde{\varepsilon}_{(g,m)}$ are defined in the same way. Further, using that $\boldsymbol{\hat{\beta}}_n = \boldsymbol{\beta}+ o_p(1)$, 
$$
\boldsymbol{\acute{\Sigma}}_n
=
\frac{M}{M-1}
\sum_{g=1}^{G_n} \sum_{m=1}^M 
\boldsymbol{\tilde{x}}_{(g,m)} \, \boldsymbol{\tilde{x}}_{(g,m)}^\prime 
\
\left( 
\varepsilon_{(g,m)} \tilde{\varepsilon}_{(g,m)} 
+
\left(\boldsymbol{x}_{(g,m)}^\prime \boldsymbol{\beta}+\alpha_g\right) \, \tilde{\varepsilon}_{(g,m)} 
\right) + o_p(G_n).
$$
Because $\mathbb{E}(\varepsilon_{(g,m)} \tilde{\varepsilon}_{(g,m)} \vert \boldsymbol{x}_{(g,1)},\ldots,\boldsymbol{x}_{(g,M)}) = \sigma^2_{(g,m)} \, (M-1)/M$ the first term constitutes an unbiased estimator of $\boldsymbol{\Sigma}_n$. The second term on the right-hand side is mean zero because the errors are mean-independent of the regressors. Its variance, however, depends on the $\alpha_1,\ldots, \alpha_{G_n}$. 


In a two-wave panel,
$$
\boldsymbol{\acute{\Sigma}}_n
=
\frac{1}{4}
\sum_{g=1}^{G_n} \Delta \boldsymbol{x}_{g} \, \Delta \boldsymbol{x}_{g}^\prime \, (\Delta y_{g} - \Delta \boldsymbol{x}_{g}^\prime \boldsymbol{\hat{\beta}}_n) \, \Delta y_g,
$$
where $\Delta$ denotes the first-difference operator; so, e.g., $\Delta y_{g} := y_{(g,2)}-y_{(g,1)}$. In this case the variance estimator does not depend on the fixed effects. The factor $\frac{1}{4}$ appears because of the de-meaning and is inconsequential. The implied estimator of the covariance matrix  $\boldsymbol{\Omega}_n$ is
$$
\left(\sum_{g=1}^{G_n} \Delta \boldsymbol{x}_g \, \Delta \boldsymbol{x}_g^\prime \right)^{-1}
\left(
\sum_{g=1}^{G_n} \Delta \boldsymbol{x}_{g} \, \Delta \boldsymbol{x}_{g}^\prime \, (\Delta y_{g} - \Delta \boldsymbol{x}_{g}^\prime \boldsymbol{\hat{\beta}}_n) \, \Delta y_g
 \right)
\left(\sum_{g=1}^{G_n} \Delta \boldsymbol{x}_g \, \Delta \boldsymbol{x}_g^\prime \right)^{-1}.
$$
The same estimator would be obtained if our procedure would be applied directly to the first-differenced model
$
\Delta y_g = \Delta \boldsymbol{x}_g^\prime\boldsymbol{\beta} + \Delta \varepsilon_g.
$
The standard estimator of \cite{Eicker1963,Eicker1967} and \cite{White1980} applied to this model is
$$
\left(\sum_{g=1}^{G_n} \Delta \boldsymbol{x}_g \, \Delta \boldsymbol{x}_g^\prime \right)^{-1}
\left(
\sum_{g=1}^{G_n} \Delta \boldsymbol{x}_{g} \, \Delta \boldsymbol{x}_{g}^\prime \, (\Delta y_{g} - \Delta \boldsymbol{x}_{g}^\prime \boldsymbol{\hat{\beta}}_n)^2
 \right)
\left(\sum_{g=1}^{G_n} \Delta \boldsymbol{x}_g \, \Delta \boldsymbol{x}_g^\prime \right)^{-1},
$$
and is known to be consistent here.

In the one-way panel model our Assumption \ref{ass:inference} holds provided that
$$
\frac{\sum_{g=1}^{G_n} \lVert \alpha_g\rVert^{2+\theta}}{G_n} = O(1)
$$
and $\max_g \max_m\mathbb{E}(\lVert\boldsymbol{x}_{(g,m)} \rVert^{2+\theta})=O(1)$ for some $\theta>0$. Given that \cite{CattaneoJanssonNewey2018} impose that $\max_g \max_m\mathbb{E}(\lVert\boldsymbol{x}_{(g,m)} \rVert^{4})=O(1)$ to validate Assumption \ref{ass:design}, the condition on the fixed effects is the only additional requirement needed for our variance estimator to be consistent in this model.

\paragraph{Linear model with increasing dimension.}
Finally, consider the regression model that takes \eqref{eq:model} as the data generating process for independent and identically distributed observations $(y_{i,n},\boldsymbol{x}_{i,n}, \boldsymbol{w}_{i,n})$, i.e.,
$$
y_{i,n} = \boldsymbol{x}_{i,n}^\prime \boldsymbol{\beta} + \boldsymbol{w}_{i,n}^\prime \boldsymbol{\gamma}_n + u_{i,n},
\qquad
n\, \mathbb{E}(\lVert \mathbb{E}(u_{i,n} \vert \boldsymbol{x}_{i,n},\boldsymbol{w}_{i,n}) \rVert^2) 
= n \varrho_n = o(1),
$$
as in \cite{CattaneoJanssonNewey2018}. The generic nature of this model makes it difficult to specify a single set of simple primitive conditions for our Assumption \ref{ass:inference} to hold. First consider the requirement that $\sum_{i=1}^n \lVert\boldsymbol{\tilde Q}_{i,n}\rVert^4 = O_p(n)$. A sufficient condition here is that
$$
n\, \chi_n 
=
n\, 
\min_{\boldsymbol{\delta}}
\mathbb{E}\left(
\lVert
\mathbb{E}(\boldsymbol{x}_{i,n} \vert \boldsymbol{w}_{i,n})
-
\boldsymbol{\delta}^\prime \boldsymbol{w}_{i,n}
\rVert^2\right) = O(1).
$$
This will be the case, for example, when $\boldsymbol{w}_{i,n}$ is discrete and a saturated regression model is used. It will also hold under smoothness conditions on the function $\mathbb{E}(\boldsymbol{x}_{i,n} \vert \boldsymbol{w}_{i,n})$ when the $\boldsymbol{w}_{i,n}$ are approximating functions, as discussed above. Alternatively, if $\chi_n = O(1)$, a sparsity condition on the projection matrix $\boldsymbol{M}_n$ can be used. One such condition is $\max_i n_{i,n}=O_p(1)$, together with $\sum_{i=1}^n \mathbb{E}(\lVert \boldsymbol{Q}_{i,n}\rVert^4) = O(n)$, where we have introduced $n_{i,n}:= \sum_{j=1}^n  \lbrace (\boldsymbol{M}_n)_{i,j} \neq 0\rbrace $.  \cite{CattaneoJanssonNewey2018} showed that the weaker condition $\max_i n_{i,n}=o_p(n^{1/3})$ can be used to support Assumption \ref{ass:approximations} under the same moment condition on the $\boldsymbol{Q}_{i,n}$. Such a rate would appear difficult to obtain here. Finally, for $\max_i \lVert \mu_{i,n} \rVert /\sqrt{n}=o_p(1)$ to hold we again impose that $\mathbb{E}(\lVert \boldsymbol{x}_{i,n} \rVert^{2+\theta}) = O(1)$ for some $\theta>0$. When $\boldsymbol{w}_{i,n}$ are approximating functions, a similar moment condition on the function that is being approximated will again suffice. In the case where $\boldsymbol{w}_{i,n}$ are just many control variables included in the regression we can again use a sparsity condition. Let $\kappa_{i,n}$ denote the number of nuisance parameters on which $y_{i,n}$ depends. Then one alternative sufficient condition is that
$$
\max_i \kappa_{i,n} = O(n^{\frac{1}{2}\frac{\theta}{2+\theta}}),
$$
for some $\theta>0$, together with the assumption that the entries of $\boldsymbol{w}_{i,n}$ have $2+\theta$ moments. A condition on $\kappa_{i,n}$ is different from a condition on $n_{i,n}$, as it only pertains to the regression of $y_{i,n}$ on $\boldsymbol{x}_{i,n}$ and $\boldsymbol{w}_{i,n}$ and does not restrict the auxiliary regression of $\boldsymbol{x}_{i,n}$ on $\boldsymbol{w}_{i,n}$. 
When the covariates are normally distributed or have bounded support, for example, we can allow for $\max_i \kappa_{i,n} = O(\sqrt{n})$.

\section{Numerical illustrations} \label{sec:4}

\subsection{Simulations}
We present numerical results for a setup taken from \cite{CattaneoJanssonNewey2018}. Data are generated as
$$
y_{i,n} = x_{i,n} \beta + \boldsymbol{w}_{i,n}^\prime \boldsymbol{\gamma}_n + \varepsilon_{i,n},
$$
where $x_{i,n}\sim \mathrm{i.i.d.}~\boldsymbol{N}(0,1)$, $\boldsymbol{w}_{i,n}$ contains a constant term and a collection of $q_n-1$ zero/one dummy variables, and $\varepsilon_{i,n}\sim \mathrm{i.i.d.}~\boldsymbol{N}(0,1)$. The dummy variables are drawn independently with success probability $\pi$ and $\boldsymbol{\gamma}_n = \boldsymbol{0}$. The sample size was fixed to $n=700$ throughout and we considered $q_n \in \lbrace 1, 71, 141, 211, 281, 351, 421, 491, 561, 631 \rbrace$. All statistics reported below were computed over $10,000$ Monte Carlo replications and all the variables were redrawn in each replication.

We consider three designs that vary in $\beta$ and $\pi$. Design A is the design of \cite{CattaneoJanssonNewey2018}. It has $\beta=1$ and $\pi=.02$.\footnote{The description of the simulation design in \citet[p.~1358]{CattaneoJanssonNewey2018} contains a typo that would imply that $\pi=.0062$.} Each dummy variable takes on the value one for about $14$ out of $700$ observations, on average. Design B is a more sparse design where $\beta=1$ is maintained and $\pi$ is reduced to $.01$, leading to each dummy switching on for only 7 observations in each replication, on average. Design C, in turn, sets $\beta=2$ and maintains $\pi=.02$. This implies that the conditional variance of $y_{i,n}\, \acute{u}_{i,n}$ increases by a factor of four.

The results for the three design variations are presented in Tables \ref{table:MC1extended}--\ref{table:MC3extended}. Empirical size of the two-sided $t$-test of the null that $\beta=1$ and the average width of the corresponding confidence interval are given for all variance estimators discussed;  HC0, its modifications HC1, HC2, and HC3, the bias-corrected estimator of \cite{CattaneoJanssonNewey2018} (HCK) and the estimator presented here (HCA). 

As the setup features homoskedastic errors, inference based on HC0 will be  liberal for large $n$ when $q_n/n$ is not small (\citealt{ChesherJewitt1987}, \citealt{CattaneoJanssonNewey2018}). This is apparent from inspection of the tables. The degrees-of-freedom correction performed by HC1 alleviates most of this concern here. 
HC2, which is consistent under homoskedasticity, performs quite similarly to HC1. Both corrections do come with (on average) wider confidence intervals. The simulations also illustrate that HC3 yields conservative inference. The rejection frequency (under the null) approaches zero as $q_n/n$ increases. The confidence intervals are also wide. This implies that tests based on HC3 will suffer from low power. The relative inefficiency compared to HC1 and HC2 grows as $q_n/n$ increases.

In Design A HCK gives close to correct size and confidence intervals of a comparable length as HC1 and HC2 for most of the values of $q_n/n$. As $q_n/n$ increases above $\frac{1}{2}$ this variance estimator does not always exist in each replication. When this happens in our Monte Carlo, HCK defaults to HC0. Consequently, for the larger values of $q_n/n$ in Table \ref{table:MC1extended}, the size of the HCK method increases above its nominal size and the average length of the confidence interval equally shrinks relative to HC1 and HC2. In the sparser Design B the non-existence of the HCK estimator is more frequent and also arises for smaller values of $q_n/n$. This explains the large overrejection rates observed in Table \ref{table:MC2extended}. The results for Design C in Table \ref{table:MC3extended} are similar to those in Table \ref{table:MC1extended} as the performance of HCK is invariant to the scale of the regression slopes.

The simulation results also confirm the theory behind our variance estimator. It yields close to correct inference in all three designs and for all values of $q_n/n$ considered, although some more overrejection is observed for the top-end values of this ratio. The average length of its confidence interval is comparable to those obtained for HC1, HC2, and HCK and is substantially smaller than those for HC3 for large values of $q_n/n$. These results show that this estimator is useful in problems with many covariates when HCK is not available, for example due to sparseness of the regressor design or because of the presence of high-leverage observations more generally.

The supplemental appendix contains additional simulation results for a partially-linear regression model and for a one-way panel model.


\subsection{Empirical example}
We next use the different variance estimators available to infer the union membership premium. The data are a balanced panel on 545 working individuals and span 8 years (1980--1987), giving a total of 4,360 observations. They are taken from \cite{VellaVerbeek1998} and are available from the data archive of the {\em Journal of Applied Econometrics} (http://qed.econ.queensu.ca/jae/1998-v13.2/vella-verbeek/).

We estimate the union premium as the coefficient from a least-squares regression of log wages on a dummy for union membership, after partialling-out a set of control variables. This set contains average hours worked per week, dummies for marital status and for poor health, a quadratic term in years of experience, and a large collection of dummies, as follows. First, as the data are a panel, both individual fixed effects and year fixed effects are included. Second, as the data contain information on the type of occupation (out of 9 categories) and the industry (out of a total of 12) in which the job is located we also control for these by including sector and industry dummies. Third, we allow for interaction effects between these categorical variables by including occupation-by-year, sector-by-year, and occupation-by-sector dummies as well as occupation-by-sector-by-year dummies. Baseline categories for the year-, sector-, and occupational dummies are chosen and their corresponding dummies are dropped as to avoid a dummy-variable trap. Certain interactions of the occupation and sector dummies never take on the value one and so are equally removed from the analysis. This gives a total of 1,086 control variables relative to 4,360 observations. Our point estimate of the union premium is $7.43\%$. This is in line with the literature (see, e.g., \citealt{Jakubson1991}).

The standard errors on this point estimate, obtained by the various methods discussed, together with the implied $95\%$ confidence intervals for the union premium are collected in Table \ref{tab:example}. A total of 330 observations have leverage that exceeds the one-half threshold. 36 of these observations can be perfectly explained by the control variables because of the inclusion of the occupation-by-sector and the occupation-by-sector-by-year dummies and so contain no information on the union premium. The HCK variance estimator could not be computed here. (We note that the invertibility condition fails also after dropping the 36 non-informative observations.) This explains why an entry for this estimator is not available. HC0 yields the smallest standard error on our point estimate ($1.72\%$). HC3 yields the largest ($2.32\%$). The standard error of HCA ($1.93\%$) is roughly in the middle of these two bounds. HC1 and HC2 give similar, slightly larger, standard errors as HCA ($1.99\%$ and $1.98\%$).


\section{Conclusion}
This paper has presented a heteroskedasticity-robust covariance-matrix estimator for linear regression models that is consistent under an asymptotic scheme where the number of control variables, $q_n$, grows at the same rate as the sample size, $n$. The estimator is similar to the proposal of \cite{KlineSaggioSoelvsten2018} but our consistency result covers more general settings. The estimator complements work by \cite{CattaneoJanssonNewey2018}, who derived inconsistency results for members of the HC-class of variance estimators, proved asymptotic-conservativeness of the HC3 estimator, and presented an alternative variance estimator (that is based on \citealt{HartleyRaoKiefer1969}) that remains consistent when $\lim\sup_n q_n/n < \frac{1}{2}$. Under a set of additional high-level conditions, our estimator allows to weaken this restriction to $\lim\sup_n q_n/n < 1$. Primitive conditions for these where given for partially-linear models, fixed-effect panel data models, and generic regression models with increasing dimension. Simulation verify the theoretical properties. An empirical application to estimation of the union premium from panel data was also presented.

The idea underlying our variance estimator can be useful as a device to correct for bias more generally. \cite{KlineSaggioSoelvsten2018} use it to bias-correct quadratic forms in fixed-effect estimators. \cite{Chernozhukov2018} and \cite{NeweyRobins2018} use related cross-fitting techniques to reduce bias in high-dimensional estimation problems that feature machine-learning estimators. \cite{CattaneoJanssonMa2018} characterize the bias in (nonlinear) two-step estimators when the first step features a high-dimensional linear regression. 


\setlength{\bibsep}{1pt} 
\bibliographystyle{chicago3}
\bibliography{hcstd}

\newpage

\begin{table}[htbp]
  \centering 
  \caption{Design A ($\beta=1$ and $\pi=.02$)}
    \begin{tabular}{lrrrrrrrrrr}
		\hline\hline
$q_n$   & 1     & 71    & 141   & 211   & 281   & 351   & 421   & 491   & 561   & 631 \\
  $q_n/n$    & .0014 &	.1014 &	.2014 &	.3014 &	.4014 &	.5014 &	.6014 &	.7014 &	.8014 &	.9014 \\ 
\multicolumn{11}{c}{Rejection frequency of $5\%$-level test} \\
HC0
& .0502 & .0671 & .0771 & .1050 & .1287 & .1605 & .2054 & .2746 & .3691 & .5309 \\
HC1
& .0495 & .0518 & .0485 & .0534 & .0531 & .0459 & .0426 & .0476 & .0438 & .0446 \\
HC2
& .0498 & .0529 & .0517 & .0567 & .0574 & .0505 & .0495 & .0556 & .0530 & .0589 \\
HC3
& .0495 & .0416 & .0280 & .0241 & .0131 & .0058 & .0022 & .0005 & .0001 & .0000 \\
HCK
& .0498 & .0530 & .0520 & .0588 & .0593 & .0563 & .0598 & .0770 & .0976 & .2688 \\
HCA
& .0519 & .0535 & .0529 & .0570 & .0566 & .0524 & .0521 & .0605 & .0567 & .0674 \\
\multicolumn{11}{c}{Average width of $95\%$ confidence interval} \\ 
HC0
& .1477 & .1484 & .1490 & .1494 & .1502 & .1510 & .1517 & .1525 & .1533  & .1541 \\
HC1
& .1479 & .1567 & .1669 & .1789 & .1944 & .2141 & .2408 & .2797 & .3452 & .4944 \\
HC2
& .1478 & .1560 & .1654 & .1766 & .1909 & .2093 & .2340 & .2698 & .3300 & .4646 \\
HC3
& .1479 & .1645 & .1851 & .2113 & .2468 & .2964 & .3707 & .4937 & .7405 & 1.4803 \\
HCK
& .1478 & .1559 & .1651 & .1758 & .1894 & .2065 & .2292 & .2606 & .3064 &  .3482\\
HCA
& .1477 & .1557 & .1652 & .1763 & .1907 & .2087 & .2334 & .2688 & .3284 & .4611 \\

\hline\hline
    \end{tabular}
  \label{table:MC1extended}
\end{table}%

\begin{table}[htbp]
  \centering 
  \caption{Design B ($\beta=1$ and $\pi=.01$)}
    \begin{tabular}{lrrrrrrrrrr}
		\hline\hline
$q_n$   & 1     & 71    & 141   & 211   & 281   & 351   & 421   & 491   & 561   & 631 \\
  $q_n/n$    & .0014 &	.1014 &	.2014 &	.3014 &	.4014 &	.5014 &	.6014 &	.7014 &	.8014 &	.9014 \\ 
\multicolumn{11}{c}{Rejection frequency of $5\%$-level test} \\
HC0
& .0479  & .0606 & .0767 & .0963 & .1137 & .1457 & .1907 & .2450 & .3353 & .4639 \\
HC1
& .0474  & .0462 & .0479 & .0471 & .0410 & .0432 & .0390 & .0340 & .0291 & .0256 \\
HC2
& .0477  & .0477 & .0526 & .0523 & .0496 & .0543 & .0553 & .0544 & .0555 & .0584 \\
HC3
& .0474  & .0376 & .0301 & .0227 & .0121 & .0050 & .0026 & .0005 & .0000 & .0000 \\
HCK
& .0477  & .0582 & .0765 & .0963 & .1136 & .1457 & .1907 & .2450 & .3353 & .4639 \\
HCA
& .0492  & .0512 & .0529 & .0540 & .0523 & .0566 & .0571 & .0553 & .0602 & .0705 \\
\multicolumn{11}{c}{Average width of $95\%$ confidence interval} \\ 
HC0
& .1479 & .1492 & .1504 & .1520 & .1537 & .1561 & .1588 & .1621 & .1671 & .1768 \\
HC1
& .1481 & .1575 & .1685 & .1820 & .1989 & .2212 & .2518 & .2971 & .3757 & .5648 \\
HC2
& .1480 & .1560 & .1653 & .1765 & .1905 & .2089 & .2335 & .2691 & .3291 & .4633 \\
HC3
& .1481 & .1645 & .1848 & .2109 & .2459 & .2953 & .3688 & .4908 & .7352 & 1.4661 \\
HCK
& .1480 & .1501 & .1507 & .1520 & .1537 & .1561 & .1588 & .1621 & .1671 & .1768 \\
HCA
& .1479 & .1557 & .1649 & .1763 & .1903 & .2085 & .2333 & .2686 & .3278 & .4590 \\

\hline\hline
    \end{tabular}
  \label{table:MC2extended}
\end{table}

\begin{table}[htbp]
  \centering 
  \caption{Design C ($\beta=2$ and $\pi=.02$)}
    \begin{tabular}{lrrrrrrrrrr}
		\hline\hline
$q_n$   & 1     & 71    & 141   & 211   & 281   & 351   & 421   & 491   & 561   & 631 \\
  $q_n/n$    & .0014 &	.1014 &	.2014 &	.3014 &	.4014 &	.5014 &	.6014 &	.7014 &	.8014 &	.9014 \\ 
\multicolumn{11}{c}{Rejection frequency of $5\%$-level test} \\
HC0
& .0505 & .0667 & .0800 & .1052 & .1273 & .1611 & .2045 & .2653 & .3708 & .5273 \\
HC1
& .0501 & .0535 & .0495 & .0508 & .0490 & .0466 & .0464 & .0488 & .0462 & .0453 \\
HC2
& .0504 & .0542 & .0519 & .0536 & .0529 & .0518 & .0521 & .0563 & .0545 & .0583 \\
HC3
& .0501 & .0420 & .0310 & .0209 & .0117 & .0062 & .0019 & .0006 & .0000 & .0000 \\
HCK
& .0504 & .0545 & .0525 & .0557 & .0553 & .0577 & .0631 & .0763 & .1026 & .2635 \\
HCA
& .0555 & .0580 & .0584 & .0612 & .0588 & .0567 & .0624 & .0673 & .0710 & .0935 \\
\multicolumn{11}{c}{Average width of $95\%$ confidence interval} \\ 
HC0
& .1480 & .1483 & .1489 & .1496 & .1502 & .1508 & .1515 & .1523 & .1535 & .1543 \\
HC1
& .1482 & .1566 & .1668 & .1792 & .1944 & .2139 & .2405 & .2794 & .3458 & .4950 \\
HC2
& .1481 & .1559 & .1654 & .1769 & .1909 & .2091 & .2336 & .2696 & .3304 & .4652 \\
HC3
& .1482 & .1644 & .1851 & .2116 & .2468 & .2961 & .3698 & .4933 & .7411 & 1.4817 \\
HCK
& .1481 & .1558 & .1650 & .1761 & .1894 & .2064 & .2286 & .2601 & .3066 & .3491 \\
HCA
& .1472 & .1552 & .1644 & .1761  & .1900 & .2080 & .2320 & .2673 & .3266 & .4535 \\

\hline\hline
    \end{tabular}
  \label{table:MC3extended}
\end{table}

\newpage

\begin{table}[htbp]
  \centering
  \caption{Inference on the union premium }
    \begin{tabular}{lcc}
		\hline\hline
	\multicolumn{3}{c}{$\hat{\beta}_n=.0743$}	
		\\
& {std.~error} & $95\%$ conf.~int. \\
    HC0   & .0172 & [.0406, .1081] \\
    HC1   & .0199 & [.0354, .1133] \\
    HC2   & .0198 & [.0356, .1131] \\
    HC3   & .0232 & [.0288, .1199] \\
    HCA   & .0193 & [.0365, .1122] \\
	\hline\hline
    \end{tabular}
  \label{tab:example}
\end{table}

\end{document}